\documentclass[letterpaper, 10 pt]{IEEEtran}
\IEEEoverridecommandlockouts

\usepackage{color}
\usepackage{enumerate,graphicx,epstopdf}
\usepackage{amsmath,amssymb,bm}
\usepackage{cite}
\usepackage{mathtools}
\usepackage{array}
\usepackage{amsthm}
\usepackage{algpseudocode,algorithm,algorithmicx}

\def\tto{\;{\lower 1pt \hbox{$\rightarrow$}}\kern -12pt
           \hbox{\raise 2.8pt \hbox{$\rightarrow$}}\;}
\def\bbs{\bar s}
\def\downto{{\raise 1pt \hbox{$\scriptstyle \,\searrow\,$}}}
\def\la{\langle}\def\ra{\rangle}
\def\bx{\bar x} \def\bq{\bar q} 
\def\bu{\bar u} 
\def\by{\bar y} \def\bv{\bar v} \def\bp{\bar p} 
  \def\bz{\bar z} \def\bp{\bar p}

\def\gph{\mathop{\rm gph}\nolimits}

\def\eps{\varepsilon}    

\newtheorem{rmk}{Remark}

\newtheorem{dfn}{Definition}

\newtheorem{theorem}{Theorem}

\newcommand{\reals}{\mathbb{R}}

\newcommand\mc[1]{\mathcal{#1}}

\title{Sensitivity-based Warmstarting for Nonlinear Model Predictive Control with Polyhedral State and Control Constraints\thanks{A. Dontchev is with the American Mathematical Society and the University of Michigan, Email: \{dontchev@umich.edu\}. D. Liao-McPherson and I. Kolmanovsky are with the University of Michigan, Ann Arbor. Email:\{dliaomcp, ilya\}@umich.edu. M. Nicotra is with the University of Colorado Boulder, Email: \{marco.nicotra@colorado.edu\}. V. Veliov is with Technical University of Vienna, Austria, Email: \{veliov@tuwien.ac.at\}.  This research is supported by the National Science Foundation Award Number CMMI 1562209, the Austrian Science Foundation (FWF) Grant  P31400-N32, and the Australian Research
Council (ARC) Project DP160100854.}}
\author{Dominic Liao-McPherson, Marco M. Nicotra, Asen L. Dontchev, Ilya V. Kolmanovsky, Vladimir. M. Veliov}

\begin{document}
\maketitle

\begin{abstract}
Model predictive control (MPC) is of increasing interest in applications for constrained control of multivariable systems. However, one of the major obstacles to its broader use is the computation time and effort required to solve a possibly non-convex optimal control problem (OCP) online. This paper introduces a sensitivity-based warmstarting strategy for systems with nonlinear dynamics and polyhedral constraints with the goal of reducing the computational footprint of MPC controllers. It predicts changes in the solution of the parameterized OCP as the parameter varies, by calculating the semiderivative of the solution mapping. The main novelty of the paper is that the polyhedrality of the constraints allows us to avoid imposing any constraint qualification conditions or strict complementarity assumptions. A numerical study featuring MPC applied to unmanned aerial vehicles illustrates the proposed approach.
\end{abstract}

\section{Introduction} \label{SIntro}

In Model Predictive Control (MPC) \cite{grune2017nonlinear,rawlings2009model} control actions are computed by solving constrained optimal control problems (OCPs) in real-time. MPC can systematically handle nonlinearities and constraints, but it requires solving (possibly approximately) a potentially non-convex OCP at each sampling instance. This motivates research into advanced numerical methods to enable MPC implementation.

At each sampling instance model predictive controllers measure or estimate the system state and then solve a discrete-time OCP, using the estimated state as an initial condition in the OCP, to determine the MPC-generated control action. As a result, the OCP is parameterized by the initial state. Often the states of the system at subsequent sampling instances are close so that, if the OCP satisfies appropriate regularity conditions, the solutions of the OCPs will be close as well. If we can determine an estimate of the change of the optimal solution, we can use this information to predict the optimal solution at the next time step and then start an optimization procedure from that prediction resulting in reduced computation time. The practice of exploiting sensitivity estimates to initialize an optimization algorithm is often referred to as sensitivity based warmstarting and is closely related to continuation/homotopy or solution tracking methods. 

Many early sensitivity methods, e.g., CGMRES\cite{ohtsuka2004continuation}, are based on continuation methods for smooth nonlinear equations and cannot directly handle inequality constraints. In \cite{zavala2009advanced} Zavala and Biegler proposed an advanced step strategy which exploits the following result established in \cite{F}: if the strong second order sufficient conditions (SSOSC), linear independence constraint qualification (LICQ) and strict complementarity slackness (SCS) condition hold, then the solution mapping is continuously differentiable in a vicinity of the solution. The derivative of the solution mapping, evaluated at the solution obtained at the previous sampling instance, can be used as a predictor for the optimal solution at the next sampling instance. However, the SCS condition is difficult to satisfy at all time instants. E.g., it cannot hold when an inequality constraint changes its activity status; typically, in such cases, the solution mapping is not differentiable with respect to the parameter. A similar method, IPA-SQP \cite{ghaemi2009integrated} computes a derivative based predictor using neighboring extremals and combines it with an sequential quadratic programming (SQP) based corrector. It handles constraint (de)activation using an active set strategy and requires that the SSOSC and LICQ hold.

In \cite{zavala2010real} Zavala and Anitescu developed a path-following strategy in the framework of parameterized generalized equations (GEs) using an augmented Lagrangian corrector-only approach. They assumed the SOSSC and LICQ but relaxed the SCS assumption, allowing the active set to vary. Similar approaches were proposed in \cite{dinh2012adjoint} using a sequential convex programming based corrector, in \cite{diehl2001real}, and in \cite{SSPC}, where a predictor and corrector are derived using tools from nonsmooth analysis. A more elaborate analysis of a path-following method for tracking solution trajectories of parameterized variational inequalities is presented in \cite{dontchev2013euler}.

If one replaces the LICQ with the weaker Mangasarian-Fromovitz constraint qualification (MFCQ) \cite{nocedal2006numerical} but requires a strengthened form of the SSOSC, i.e., that the SSOSC holds for all Lagrange multipliers, then the optimal solution can be shown to be directionally differentiable \cite{ralph1995directional}. This is exploited in a paper by J\"aschke,  Yang and  Biegler \cite{jaschke2014fast} to enhance the advanced step warmstart by relaxing the SCS assumption used in \cite{zavala2009advanced}. However, this method requires solving additional linear and quadratic programming problems to handle the non-uniqueness of the Lagrange multipliers. A survey on sensitivity and solution tracking methods can be found in \cite{wolf2016fast}.

To the best of our knowledge, all existing sensitivity based methods require a constraint qualification of some sort, e.g., \cite{zavala2009advanced,zavala2010real,ghaemi2009integrated,dinh2012adjoint} require the LICQ and \cite{jaschke2014fast} requires the MFCQ. Constraint qualifications are difficult to verify, both a-priori and a-posteriori, because they require a very accurate estimate of the solution. Moreover, the absence of a constraint qualification can lead to numerical difficulties for most optimization algorithms; in extreme cases this can lead to failure of the optimization routine and the associated MPC controller.

In this paper, we present a novel warmstarting strategy based on a sensitivity analysis of the OCP's parameter to solution mapping; our predictor is based on the Bouligand or B-derivative \cite{robinson1987local}, also known as the semiderivative; it reduces to the standard derivative when the SSOSC, LICQ, and SC all hold. This approach requires no constraint qualification, only a numerically verifiable second order sufficient condition. In exchange, we restrict ourselves to the case where the dynamics are nonlinear but the state and control constraints are represented by convex polyhedra.

\subsection{Some useful mappings} \label{ss:mappings}
This paper makes extensive use of several kinds of set valued mappings. A set-valued mapping $\mathcal{F}$ acting between $\reals^k$ and $\reals^l$ is denoted as ${\cal F}:\reals^k\tto \reals^l$, to distinguish it from a function $f:\reals^k\to\reals^l$,  while its inverse is defined as $y \mapsto {\cal F}^{-1}(y) = \{x \mid y\in {\cal F}(x)\}$. Given a closed convex set $C \subseteq \reals^n$, the tangent cone to $C$ at a point $x \in C$ is the set of all $v$ such that
$\frac{1}{\eps_k}(x^k - x) \to v $ for some $x^k \to x, x^k \in C, \eps_k \downto 0$ the normal cone mapping of $C$ is defined as
\begin{equation*}
  N_C(v) = \begin{cases}
  \{y~ |~ y^T(w-v) \leq 0 \ \  \forall w \in C\} & \text{if}~v \in C,\\
  \emptyset & \text{otherwise,}
  \end{cases}
\end{equation*}
The polar of a closed, convex cone $K$ is
\begin{equation}
K^\circ= \{y~|~\langle y,x\rangle \leq 0, \forall x \in K\},
\end{equation}
and then the tangent cone $T_C(v) = N_C^\circ(v)$. The euclidean projection onto the set $C$ is denoted by $\Pi_C(\cdot)$ and set addition/subtraction is defined as
\begin{equation}
  K_1 \pm K_2 = \{z \mid z = z_1\pm z_2,~~  z_1 \in K_1, z_2 \in K_2\}.
\end{equation}
For any $x \in C$ and $v \in N_C(x)$ the critical cone to $C$ at $x$ for $v$ is defined as
\begin{equation}
  K_C(x,v) = \{y~|~ y\in T_C(x),~ y^Tv = 0\}.
\end{equation}
Now suppose $C$ is polyhedral; then, by definition, there exists a matrix $\Gamma$ and a vector $b$ of appropriate dimensions such that
\begin{equation}
  C = \{x~|~\Gamma x \leq b\}.
\end{equation}
To obtain a computationally tractable expression for inclusions of the type $y\in K_C(x,v)$, define the active constraint set,
\begin{equation}
  \mathcal{A}(x) = \{i\in [1,l]~|~ \Gamma_i x = b_i\},
\end{equation}
where $l$ is the number of rows in $\Gamma$. Then the critical cone can be expressed as
\begin{equation}\label{eq:CriticalCone}
   \begin{array}{ll}
   K_C(x,v) = \{y~|&\Gamma_i y \leq 0,~i \in \mathcal{A}(x), ~y^T v = 0\},
   \end{array}
\end{equation}
see e.g., \cite[Theorem 2E.3]{book}. This can be further simplified by noting that a constraint can be deactivated if it is locally redundant with respect to the other constraints. Given an active constraint $i \in \mathcal{A}(x)$ define the polyhedral set 
\begin{equation} \label{eq:Ci_def}
C_i= \{y~|~\Gamma_j y \leq b_j,~j\in[1,l]\setminus i\}
\end{equation}
as the set that satisfies all other constraints. Due to the constraint $y^T v = 0$ in \eqref{eq:CriticalCone} if $v+ N_{C_i}(x)\ni0$, the constraint $i$ is redundant, meaning that the critical cone will remain unchanged if one were to ignore the constraint $\Gamma_i y \leq 0$. By defining the set of redundant constraints
\begin{equation}
  \bar{\mathcal{A}}(x,v) = \{i\in \mathcal{A}(x)~|~ v + N_{C_i}(x)\ni0\},
\end{equation}
the critical cone \eqref{eq:CriticalCone} can be rewritten as
\begin{equation}
   \begin{array}{ll}
   K_C(x,v)  = \{y~|&\Gamma_i y \leq 0,~i \in \bar{\mathcal{A}}(x,v),~y^T v = 0\\
   &\Gamma_i y = 0,~i \in \mathcal{A}(x)\setminus\bar{\mathcal{A}}(x,v)\},
   \end{array}
\end{equation}
where all the non-redundant constraints are now treated as equalities. In practice, the set $\bar{\mathcal{A}}(x,v)$ can be obtained by checking if $v =\Pi_{C_i}(v)$ for each $i\in\mathcal{A}(x)$ using the relationship between normal cones are projections \cite{book}. This procedure is summarized in Algorithm~\ref{algo:AS}.

\begin{algorithm}[H] 
\caption{Determine redundant constraints} \label{algo:AS}
\begin{algorithmic}[1] 
 \renewcommand{\algorithmicrequire}{\textbf{Input:}}
 \renewcommand{\algorithmicensure}{\textbf{Output:}}
 \Require $x,v$
 \Ensure  $\bar{\mathcal{A}}$
 \State $\bar{\mathcal{A}} = \emptyset$
 \ForAll{$i \in \mathcal{A}(x)$}
    \If{$v = \Pi_{C_i}(v)$}
        \State $\bar{\mathcal{A}} \gets \bar{\mathcal{A}} \cup \{i\}$ \Comment{$C_i$ is defined in \eqref{eq:Ci_def}}
    \EndIf
 \EndFor
 \end{algorithmic}
 \end{algorithm}

\section{Problem formulation}
We consider  the following discrete-time OCP:
\begin{subequations} \label{oc} 
\begin{equation}
\underset{x,u}{\min} \  \quad J(x,u) = \phi(x_N) + \sum_{i=0}^{N-1} \ell(x_i,u_i) 
\end{equation}
subject to
\begin{eqnarray}
 & x_{i+1} =  f(x_i,u_i),~~ i = 0,\dots,N-1, ~~ x_0  =  p, \label{eq:eq_constr}
  \\
&(x_i,u_i) \in Z_i,~ i = 1,\dots,N-1,\\
&u_0 \in U_0, ~ x_N \in X_N
\end{eqnarray}
\end{subequations}
where $N$ is a natural number denoting the discrete-time horizon, $x_i \in \reals^n,~u_i \in \reals^m$, $x = (x_1,\dots, x_N)$ and $u = (u_0,\dots, u_{N-1})$ are the discrete-time state and control sequences, and the initial state is regarded as a parameter $p \in \reals^n$. The functions  $\ell: \reals^n\times\reals^m \to \reals, $ $\phi:  \reals^n\to \reals$, and $f: \reals^n\times\reals^m \to \reals^n$  are  assumed  twice continuously differentiable everywhere for simplicity. The families of closed and convex sets $Z_i \subseteq \reals^{n+m}$, $U_0 \subseteq \reals^m$ and $X_N \subseteq \reals^n$ describe state and control constraints that may vary in time. Throughout the paper we assume that each of these sets is a {\em polyhedral set}; that is, it can be described by linear inequality constraints of the form
\begin{subequations} \label{eq:polydef}
\begin{gather} 
  Z_i = \{(x,u)~|~ E_i [x_i^T~~u_i^T]^T \leq c_i\},\\
  X_N = \{x~|~E_N x \leq c_N\},~~ U_0 = \{u~|~E_0 u \leq c_0\}
\end{gather}
\end{subequations}
for some matrices $E_i$ and vectors $c_i$ of compatible dimensions. We also assume that for a fixed reference value $\bp$ of the parameter problem \eqref{oc} has a solution $(\bx, \bu)$.

\begin{rmk} In this paper we consider only the initial state as a parameter for brevity. The sensitivity analysis presented can be easily extended to OCPs in which the functions in the cost and state equations also depend on additional parameters, e.g., a target state or a previewed time-varying input, if the problem functions are continously differentiable with respect to the parameters. Since these parameters often enter the OCP nonlinearly, e.g., through a quadratic term in the cost function, our approach is more general than the generalized tangential predictor described in Section 5.3 of \cite{magni2009nonlinear} which requires that the parameter enter linearly.
\end{rmk}

\section{Optimality and sensitivity }\label{back}

Problem \eqref{oc} is a special case of the following  parameterized optimization problem
\begin{subequations} \label{opt}
\begin{gather} 
  \underset{v}{\min}\quad  J(v)\\
   \text{ subject to } g(p,v) = 0, \quad v \in V,
\end{gather}
\end{subequations}
where  $v= (u_0,x_1,u_1~\ldots,~u_{N-1},x_N) \in \reals^j,$ $ j=N(m+n)$, $J:\reals^j \to\reals$ is defined as in \eqref{oc}, the function $g: \reals^{n+j} \to \reals^d$, $d= nN$, is given by
\begin{equation} 
  g(p,v) = \begin{bmatrix}
  x_1 - f(p,u_0)\\
  x_2 - f(x_1,u_1)\\
  \vdots\\
  x_N - f(x_{N-1},u_{N-1})
  \end{bmatrix},
\end{equation}
and where $V = U_0 \times Z_1 \times \hdots \times Z_{N-1}\times X_N$. The associated  Lagrangian has the form
\begin{equation}
  \mc{L}(p,v,q) = J(v) + q^T g(p,v),
\end{equation}
where $q\in \reals^{Nn}$ is the vector of  the Lagrange multipliers associated with the equality constraints in \eqref{opt}; in the context of optimal control it is usually called the vector of costates.

Let $\bv = (\bu_0,\bx_1,\bu_1,\ldots, \bu_{N-1},\bx_N)$ be a local minimizer of \eqref{opt} for a reference value $\bp$ of the parameter. It is known, see e.g., \cite[Theorem 6.14]{rockafellar2009variational} and \cite[Theorem 5.1.1]{ioffe2009theory}, that under the constraint qualification condition
\begin{equation} \label{cc}
\text{the matrix}  \nabla_v g(\bar p, \bar v) \in \reals^{Nn \times N(m+n)}   \text{ is surjective} 
\end{equation}
the first-order necessary conditions for optimality are 
\begin{gather} \label{kkt}
\nabla_v \mathcal{L}(p,v,q) + N_V(v) \ni 0,\\
g(p,v) = 0.
\end{gather}
Noting that $\nabla_q {\cal L} = g$, defining $z = (v,q) $, $E= V\times \reals^d$, and letting $F(p,z) = \nabla_z {\cal L}(p,v,q)$,
we arrive at the following parameterized variational inequality (VI):
\begin{equation} \label{vi}
F(p,z) + N_E(z) \ni 0.
\end{equation}
It is easy to show that \eqref{cc} holds for problem \eqref{oc}. Indeed, denoting
 $B_0 = \nabla_{u_0} f(\bp, \bu_0),$
$A_i = \nabla_{x_i} f(\bx_i,\bu_i)$, and $B_i = \nabla_{u_i} f(\bx_i,\bu_i), i = 1, \dots, N-1,$ 
the surjectivity of the matrix $\nabla_v g(\bp,\bv)$  becomes the condition that for every $\xi = (\xi_0, \dots,\xi_{N-1})\in \reals^{Nn}$ the system 
\begin{gather}
x_1 - B_0u_0  =  \xi_0,\\
 x_{i+1}- A_ix_i - B_iu_i  =  \xi_{i}, \quad i = 1, \dots, N-1, 
\end{gather}
has a solution. This condition clearly holds: simply choose an arbitrary sequence $(u_0, \dots, u_{N-1})$ and determine $(x_1, \dots, x_{N})$ recursively.

The solution mapping of \eqref{vi} is 
\begin{equation}
  p \mapsto S(p) = \{z~|~ F(p,z) + N_E(z) \ni 0\}.
\end{equation}
There is a well-developed theory for the properties of $S(p)$, most of which is collected in \cite[Chapter 2]{book}. We will use the following definitions:

\begin{dfn} (Strong regularity)
A set-valued mapping  ${\cal F}: \reals^k \tto \reals^l$ is said to be strongly regular at $\bx$ for $\by$ if $\by \in {\cal F}(\bx)$ and the inverse ${\cal F}^{-1}$ has a Lipschitz localization around $\by$ for $\by$; that is,  there exist neighborhoods $U$ of $\bx$ and $V$ of $\by$ such that the truncated inverse mapping $V \ni \by  \mapsto {\cal F}^{-1}(y)\cap U$ is  a function which is Lipschitz continuous around $\by$.
\end{dfn}

\begin{dfn} (Semidifferentiability)
A function $\psi:\reals^k\to\reals^l$ is said to be {\em semidifferentiable} at $\bp$ if there exists a positively
homogeneous function  $D\psi(\bp): \reals^k \to \reals^l$ with the property that  for every $\eps>0$ there exists $\delta >0$ such that
for every $\Delta p \in \reals^k$ with  $\|\Delta p\|\leq \delta$
one has
$$\|\psi(\bp+\Delta p) -\psi(\bp) - D\psi(\bp)(\Delta p)\|\leq \eps\|\Delta p\|.$$
\end{dfn}

When $D\psi(\bp)$ happens to be linear, then it becomes the usual (Fr\'echet) derivative of $\psi$ at $\bp$. A semidifferentiable function is, in particular, directionally differentiable with the derivative in the direction of $h$ satisfying $\psi'(\bp;h)= D\psi(\bp)(h)$. The opposite implication holds if the function $\psi$ is Lipschitz continuous. For more about semidifferentiable functions, including examples, see e.g. \cite[Section 2D]{book}. 

The {\it  critical cone} for \eqref{vi} at $(p, z) \in \gph S$ is 
$$ K(p,z) := K_E(z,F(p,z)) = \{ w \in T_E(z) \mid w^T F(p, z) = 0\},$$
see Section~\ref{ss:mappings} for more information on critical cones. Accordingly, the subspace defined as 
\begin{equation} \label{eq:Kplus}
  K^+(p,z) =  K(p,z) - K(p,z)
\end{equation}
is the smallest subspace that includes the critical cone, while the subspace
\begin{equation} \label{eq:Kminus}
  K^-(p,z) = K(p,z) \cap [-K(p,z)]
\end{equation}
is the largest subspace that is included in $K(p,z)$.

Using these concepts, we state the following theorem, which
is a compilation of \cite[Theorems 2E.6, 2E.8]{book}, and characterizes the strong regularity and the semidifferentiability properties of the solution  mapping of \eqref{vi}:

\begin{theorem}\label{thm1}  Let $\Lambda = \nabla_z F(\bp, \bz)$ and let $K=K_E(\bz,F(\bp,\bz))$ be the corresponding critical cone. Then suppose that the mapping $\Lambda+N_K$ is strongly regular at $0$ for $0$, this being equivalent to the condition that the linear variational inequality
\begin{equation} \label{lvi}
\Lambda z +N_K(z) \ni r
\end{equation}
has a unique solution $\bbs(r)$ for each $r \in \reals^{j+d}$. Then the solution mapping $S$ of \eqref{vi} with values $S(p)$ for $p = \bp+\Delta p$  has a Lipschitz localization $s$ at $\bp$ for $\bz$ which is
semidifferentiable and its semiderivative $Ds(\bp)(\Delta p)$ is the solution of the following variational inequality:
\begin{equation} \label{semid}
\Lambda z + N_K(z) + \nabla_p F(\bp,\bz)\Delta p \ni 0.
\end{equation}
Furthermore, in terms of the critical subspaces $K^+ = K^+(\bp, \bz)$ and $K^- = K^-(\bp, \bz)$, defined in \eqref{eq:Kplus} and \eqref{eq:Kminus}, a sufficient condition for strong regularity of $\Lambda+N_K$ or, equivalently, single-valuedness of $(\Lambda+N_K)^{-1}$, is as follows:
\begin{equation} \label{scsr}
w \in K^+, \ \Lambda w \perp K^-, \  \la w, \Lambda w\ra \leq 0 \implies w=0.
\end{equation}
\end{theorem}

Thus, in order to determine the change of the solution for a given variation $\Delta p$ of the parameter, one could compute the
corresponding semiderivative, which in turn reduced to finding the critical cone $K$ and  solving a linear VI. 

Let $\Delta p$ be a change of the parameter (the initial state). From Theorem~\ref{thm1} it follows that the semiderivative
\begin{equation}
  Ds(\bp)(\Delta p) = \begin{bmatrix}
    \Delta v^T& \Delta q^T
  \end{bmatrix}^T
\end{equation}
is a solution of the  linear VI
\begin{equation} \label{eq:lvi_ocp}
\begin{bmatrix}
  R & G^T \\
  G & 0
\end{bmatrix} \begin{bmatrix}
  \Delta v\\ \Delta q
\end{bmatrix} + N_W(\Delta v, \Delta q)  \ni \begin{bmatrix}
r_1 \\ r_2
\end{bmatrix}
\end{equation}
where $R = \nabla_{v}^2 \mc{L}(\bz,\bp)$, $G = \nabla_v g(\bv,\bp)$, $W = K_V \times \reals^n$,  where the
critical cone  $ K_V$ is
\begin{equation} \label{eq:kkcone}
   K_V = \{v \in T_V(\bv) ~ | ~  \nabla_v {\cal L}(\bp,\bz)^Tv  = 0\},
\end{equation}
and
\begin{equation} \label{eq:r_def}
  r = \begin{bmatrix}
  r_1\\
  r_2
  \end{bmatrix} = \begin{bmatrix} 
    -\nabla_{vp} \mc{L}(\bp,\bz)\\
    -\nabla_p g(\bp,\bv)
  \end{bmatrix} \Delta p.
\end{equation} 

In order to apply Theorem~\ref{thm1}, we need to adapt 
the sufficient condition for strong regularity \eqref{scsr} to our case.  Denote
$$ T= \begin{bmatrix}
  R & G^T \\
  G & 0
\end{bmatrix}$$ 
and let $W^+ = W - W$ and $W^- = W \cap [-W]$ be the critical subspaces associated  with the cone $W$.
Then condition \eqref{scsr} becomes
\begin{equation} \label{scsroc}
w \in W^+,~~ Tw \perp W^-,~~ \la w, Tw \ra \leq 0 
\implies w = 0.
\end{equation}
One should note that condition \eqref{scsroc} is a special case of 
the more elaborate {\em critical face condition} obtained in \cite{dontchev1996characterizations}, which characterizes the strong regularity  of solution mappings associated with variational inequalities over polyhedral convex sets. Obtaining a sharp, numerically tractable, form of the critical face condition for state and control constrained discrete time optimal control problems is beyond the scope of the present paper and is left for future research. However, it is possible to derive a verifiable sufficient condition, observe that \eqref{scsroc} holds provided that the matrix $T$ is positive definite on $W^+$. This in turn holds if the matrix $R$ is positive definite on the null-space of $G$; that is 
\begin{equation} \label{eq:HG_null}
  v^T R v > 0~~ \mathrm{for \ all} ~ v\neq 0 ~\mathrm{such \ that  }~~ Gv = 0.
\end{equation}
This condition is a standard second order condition in optimization\cite{nocedal2006numerical} and is always satisfied, for example, when the cost function is strongly convex and the system is linear.

The following statement summarizes the procedure for computing the semiderivatives of the solution mapping.
\begin{theorem} \label{thm:direc_deriv}
Let $\bz = (\bv,\bq)$ satisfy $\bz \in S(\bp)$, let $\Delta p$ be a variation of the parameter and assume the  condition  \eqref{eq:HG_null} holds at $(\bp,\bz)$. Then the solution mapping $p \mapsto S(p)$ has a Lipschitz localization $s(p)$ at $\bp$ for $\bz = (\bv,\bq)$ which is semidifferentiable at $\bp$ and the corresponding semiderivative $Ds(\bp)(\Delta p) = (\Delta z,\Delta q)$, or, equivalently, the directional derivative $s'(\bp;\Delta p)$, is the unique solution of the linear VI \eqref{eq:lvi_ocp}.
\end{theorem}
In Theorem~\ref{thm:direc_deriv} we assume that \eqref{eq:HG_null}, which implies \eqref{scsroc}, is satisfied at the reference solution. This condition can be enforced by appropriately regularizing the cost function \cite{bieglerIFAC}. Moreover, it's possible to monitor if \eqref{eq:HG_null} holds by checking if
\begin{equation} \label{eq:rhess_pd}
  Z^TRZ \succ 0,
\end{equation}
where the columns of $Z$ form a basis for the nullspace of $G$ and $(\cdot) \succ 0$ denotes positive definiteness. This is straightforward to check numerically by, e.g., using a QR decomposition of $G^T$ to form $Z$ and attempting to compute a Cholesky factorization of $Z^TRZ$, see \cite[Section 16.1]{nocedal2006numerical} for more details.

\begin{rmk}
The Generalized Tangential Predictor (GTP) \cite{diehl2001real}\cite[Section 5.3]{magni2009nonlinear} applied to our OCP \eqref{opt} is the directional derivative of the solution mapping of the following parameterized nonlinear program
\begin{equation} \label{eq:NLP}
\underset{v}{\mathrm{min.}}~~ J(v) ~~ \mathrm{s.t.}~~ g(v,p) = 0, ~~ Mv \leq h,
\end{equation}
where $V = \{v~|~Mv\leq h\}$. The GTP is of the form $(p,\Delta p) \mapsto (\Delta v(p,\Delta p),\Delta q(p,\Delta p),\Delta \mu(p,\Delta p))$, where $v$ are the decision variables and $q$ and $\mu$ are the dual variables associated with the equality and inequality constraints respectively. If the SSOSC and LICQ hold at a KKT point $(\bar{p},\bar{v},\bar{q},\bar{\mu})$ of \eqref{eq:NLP}, then the mapping $(\bar{p},\Delta p) \mapsto (\Delta v(\bar{p},\Delta p),\Delta q(\bar{p},\Delta p))$, i.e., the first two components of the GTP, coincides with our proposed predictor. The advantage of our predictor over the GTP is that it is well defined also in the case when the LICQ does not hold.
\end{rmk}

\begin{rmk}
The main advantage of our predictor is that it does not require any constraint qualification. Other schemes require a constraint qualification, typically the LICQ, to ensure regularity or uniqueness of the Lagrange multipliers associated with inequality constraints. Its easy to construct a situation where the LICQ does not hold by e.g., duplicating a constraint or enforcing a constraint of the form $0\leq u \leq 0$. However, the LICQ cannot be expected to always hold even in the case of linearly independent polyhedral constraints. For example, if the number of simultaneously active constraints is larger than the number of control inputs, the gradients of the of the equality constraints \eqref{eq:eq_constr} and active inequality constraints can easily become linearly dependent.
\end{rmk}

\section{A predictor-corrector algorithm} \label{ss:predictor}

Theorem~\ref{thm:direc_deriv} shows that the semiderivatives of the solution mapping can be computed by solving a linear VI. Note that \eqref{eq:lvi_ocp} are the first-order necessary conditions (i.e., the optimality system) for the following quadratic program (QP),
\begin{subequations} \label{eq:qp1}
\begin{gather}
\underset{\Delta v}{\min} \ \quad \frac12\Delta v^T R \Delta v - r_1^T\Delta v,\\
\mathrm{subject \ to}\quad G\Delta v = r_2,\\
\Delta v \in K_V, \label{eq:qp_kcone}
\end{gather}
\end{subequations}
where $R$, $G$, $K_V$, $r_1$ and $r_2$ are defined in \eqref{eq:lvi_ocp}, \eqref{eq:kkcone} and \eqref{eq:r_def}. The critical cone constraint \eqref{eq:qp_kcone} can be simplified by recalling that $K_V$ can be expressed in terms of the index set of the active constraints, see Section~\ref{ss:mappings}. In addition, recall that we can express the set $V$ as $V = \{v~|~ Mv \leq h\}$ (see \eqref{eq:polydef}) where
\begin{equation}
  M = \begin{bmatrix}
  E_0 & & \\
      & E_1 & \\
      & & \ddots &\\
      & & & E_N
  \end{bmatrix}  \text{ and } h = \begin{bmatrix}
    c_0\\c_1\\ \vdots \\c_N
  \end{bmatrix}.
\end{equation}

The QP can thus be rewritten as
\begin{subequations} \label{eq:qp3}
\begin{gather}
\underset{\Delta v}{\min} \ \quad  \frac12\Delta v^T R \Delta v + (P\Delta p)^T\Delta v,\\
\mathrm{subject \ to}\quad G\Delta v + Q\Delta p = 0, \label{eq:qp_costate}\\
M_i \Delta v \leq 0,~~ i \in \bar{\mathcal{A}},\\
M_i \Delta v = 0,~~ i \in \mathcal{A}(\bv)\setminus\bar{\mathcal{A}},\\
\nabla_v \mc{L}(\bp,\bz)^T \Delta v = 0,
\end{gather}
\end{subequations}
where $\mathcal{A}(\bv)$ and $\bar{\mathcal{A}} = \bar{\mathcal{A}}(\bv,F(\bp,\bz))$ are defined in Section~\ref{ss:mappings},  $P = \nabla_{pv}L(\bp,\bz)$, $Q = \nabla_p g(\bp,\bv)$ and $\bz = (\bv,\bq) \in S(\bp)$. The solution of \eqref{eq:lvi_ocp} can be obtained by solving \eqref{eq:qp3} and extracting the primal solution and the multipliers associated with \eqref{eq:qp_costate}. Note that \eqref{eq:qp3} requires identification of the set of active constraints but does not require any notion of strongly/weakly active constraints. Identifying a constraint as strongly or weakly active requires knowledge of a dual variable associated with that inequality. This can be problematic at points where the LICQ does not hold since then identifying the set of strongly active constraints becomes challenging,  see e.g., \cite{oberlin2006active} and the references therein.

This QP may be  difficult to solve in practice since it is not necessarily convex. This can be addressed by modifying the Hessian as follows:
\begin{equation}
  R \gets R + \rho G^TG,
\end{equation}
where $\rho \in (0,\infty)$. If \eqref{eq:HG_null} is satisfied then the modified $R$ will be positive semidefinite provided $\rho$ is chosen large enough \cite[Proposition 4.8]{izmailov2014newton}. This modification of the Hessian changes the Lagrange multipliers associated with the co-state, but the original multipliers can be recovered as $\Delta q \gets \Delta q^* - \rho Q\Delta p$ where $\Delta q^*$ are obtained by solving \eqref{eq:qp3} with $R \gets R + \rho G^TG$, and extracting the multipliers associated with \eqref{eq:qp_costate}, see e.g., \cite[Section 4.2]{izmailov2014newton}.

The predictor allows constraints that were previously active to de-activate. After the predictor step a corrector may be needed to remove constraint violations. For the corrector we use the Josephy-Newton (JN) method. Given a fixed parameter value $p$, a VI of the form \eqref{vi}, and an initial guess $z_0$, the JN method constructs an iterative sequence by repeatedly solving the following linearized VI
\begin{equation} \label{eq:JN}
  \nabla_zF(p,z_k) (z_{k+1}-z_k) + N_E(z_{k+1}) \ni 0,
\end{equation}
for $z_{k+1}$. It is well known that strong regularity of the VI implies local quadratic convergence of the JN method \cite{josephy1979newton}. 
Denoting $\delta v = v^{k+1} - v^k,$ and  $\delta q = q^{k+1} - q^k,$  \eqref{eq:JN} becomes
\begin{subequations} \label{eq:qpnc4}
\begin{multline}
\nabla_v^2\mc{L}(p,v_k,q^k)\delta v + \nabla_v\mc{L}(p,v_k,q^k) + \\\nabla_vg(p,v_k)^T \delta q + N_{V}(v_k + \delta v) \ni 0,
\end{multline}
\begin{equation}
\nabla_vg(p,v_k)\delta v + g(p,v_k) = 0.
\end{equation}
\end{subequations}
Note that this is the optimality system of the following QP:
\begin{subequations} \label{eq:qp4}
\begin{gather}
\underset{\delta v}{\mathrm{min}}\quad  \frac12\delta v^T \nabla_v^2\mc{L}(p,v_k,q_k)\delta v + \nabla_v\mc{L}(p,v_k,q_k)^T \delta v,\\
\mathrm{subject \ to}\quad \nabla_vg(p,v_k)\delta v + g(p,v_k) = 0, \label{eq:qp4_costate}\\
M \delta v \leq h - M v_k.
\end{gather}
\end{subequations}
When the Josephy-Newton iteration \eqref{eq:JN} is determined by solving the corresponding quadratic problem \eqref{eq:qp4}, this method is called Sequential Quadratic Programming (SQP). We can now summarize the predictor-corrector algorithm as Algorithm~\ref{algo:PC}. At each sampling instance the system state is measured, and the sensitivity based predictor is used to estimate the updated iterate corresponding to the parameter change based on the solution of the OCP at the previous sampling instance. This estimate is then passed to an JN based corrector loop which stops when the norm of the residual
\begin{equation}
   \pi(p,z) = ||z- \Pi_E[z - F(p,z)]||,
\end{equation}
is within the specified tolerance. The control input is then extracted from the solution and applied to the system.

Note that both the predictor and corrector steps are realized by solving QPs. The predictor QP usually has significantly fewer constraints than the corrector QP, see Section~\ref{ss:mappings}, which can lead to reduced computation times. In addition,  an initial feasible guess for the predictor QP is available, which can be helpful for both primal active-set and primal-barrier interior point methods.

\begin{rmk}
The predictor-corrector method described in Section~\ref{ss:predictor} does not introduce dual variables associated with the polyhedral inequality constraints in \eqref{opt} and thus does not require any constraint qualifications. However, the methods used to solve \eqref{eq:qp3} and \eqref{eq:qp4} for the predictor and corrector steps may introduce dual variables internally. Thus if the LICQ does not hold it could potentially cause issues for the QP solver. Primal active set methods, see e.g., \cite{nocedal2006numerical}, can encounter issues if the LICQ does not hold as can primal-dual semismooth methods \cite{liao2018regularized}. However, many implementations include heuristics or regularization schemes that handle the issue. Dual active set methods, e.g., \cite{goldfarb1983numerically} and interior point methods, see \cite{nocedal2006numerical}, typically do not require any constraint qualifications beyond Slater's condition.
\end{rmk}

\begin{algorithm}[H]
\caption{Sensitivity based Predictor-Corrector MPC} \label{algo:PC}
\begin{algorithmic}[1] 
 \renewcommand{\algorithmicrequire}{\textbf{Input:}}
 \renewcommand{\algorithmicensure}{\textbf{Output:}}
 \Require $\eps$, $p_k$, $p_{k-1}$, $v_{k-1}$, $q_{k-1}$, $\kappa$
 \Ensure  $z_k = (v_k,q_k)$
  \State Measure $p_k$
  \State $\Delta p_k = p_k - p_{k-1}$
  \State $z \gets (q_{k-1},v_{k-1})$
  \State Solve \eqref{eq:qp3} with $(\bv,\bq,\bp) = (v_{k-1},q_{k-1},p_{k-1})$ to obtain $\Delta z = (\Delta v,\Delta q)$
  \State $z \gets z + \Delta z$
  \While{$||\pi(p_k,v,q)|| > \eps$}
    \State Solve \eqref{eq:qpnc4} with $(z_k,p_k) = (z,p_k)$ to obtain $(\delta v,\delta q)$
    \State $z\gets z + \delta z$
  \EndWhile
  \State $z_k \gets z$
  \State Extract $u_0$ from $z_k$ and apply it to the plant
 \end{algorithmic}
 \end{algorithm}

\section{ Numerical simulations} \label{ss:sims}
To illustrate the theoretical  results, we consider the following model of rotational and translational dynamics of an Unmanned Aerial Vehicle (UAV),
\begin{equation}
\begin{cases}
\dot{p} = v, & m\dot{v} = T~R(\theta) e_3 - mg e_3,\\
\dot{\theta} = R(\theta)\omega, & J \dot{\omega} = - \omega \times (J \omega) + \tau,
\end{cases}
\end{equation}
where the state vector is given by the position $p\in\reals^3$, velocity $v\in\reals^3$, attitude $\theta\in(-\pi,\pi]^3$, and angular velocity $\omega\in\reals^3$ of the UAV, whereas the input vector is composed of the total thrust $T\in\reals^+$ and torques $\tau\in\reals^3$ generated by the propellers. The mass and inertia matrix of the UAV are $m=2~kg$ and $J=\text{diag}([0.82~0.82~1.62])\times10^{-2}~kgm^2$, respectively, $e_3=[0~0~1]^T$, and
\begin{equation}
  R(\theta) = \begin{bmatrix}
    1 & \sin(\theta_3) \tan(\theta_2) & \cos(\theta_3) \sin(\theta_2)\\
    0 & \cos(\theta_3) & -\sin(\theta_3) \\
    0 & \sin(\theta_3)/\cos(\theta_2) & \sin(\theta_3)/\cos(\theta_2)
  \end{bmatrix}
\end{equation}
is the attitude kinematic matrix. The system dynamics are discretized using the forward Euler approximation with $T_s=0.075~s$. The MPC is designed using the quadratic cost function,
\begin{equation}
  J(x,u)=\frac12 x_N^TPx_N + \sum_{i=0}^{N-1} \frac12 x_i^TQx_i+\frac12 u_i^TRu_i,
\end{equation}
where $Q=\text{diag}([5~5~5~10~10~10~0.1~0.1~0.1~1~1~1])$, $R=\text{diag}([0.1~0.01~0.01~0.01])$, and $P$ was computed from these values using the LQR terminal cost for the system linearized about the origin. The system is subject to the input constraints $T\in[18,22]$, $|\tau_j|\leq0.06$, as well as the state constraints $|v_j|\leq2$, where $j\in\{1,2,3\}$ refers to each component of the vector. 

\begin{figure}
  \centering
  \includegraphics[width=0.95\columnwidth]{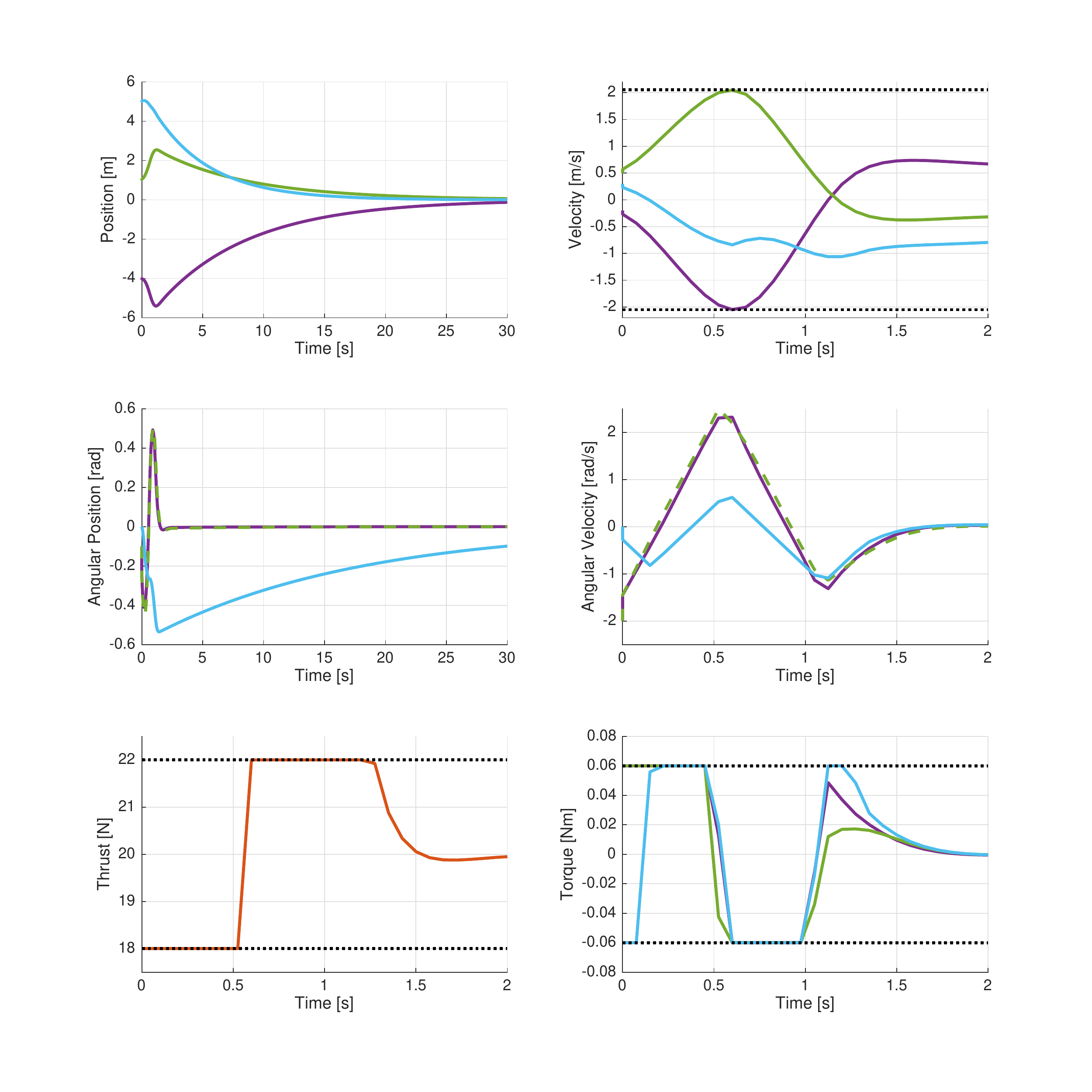}
  \caption{The closed-loop response of the UAV. All state and input constraints are respected.}
  \label{fig:response}
\end{figure}

\begin{figure}
  \centering
  \includegraphics[width=0.95\columnwidth]{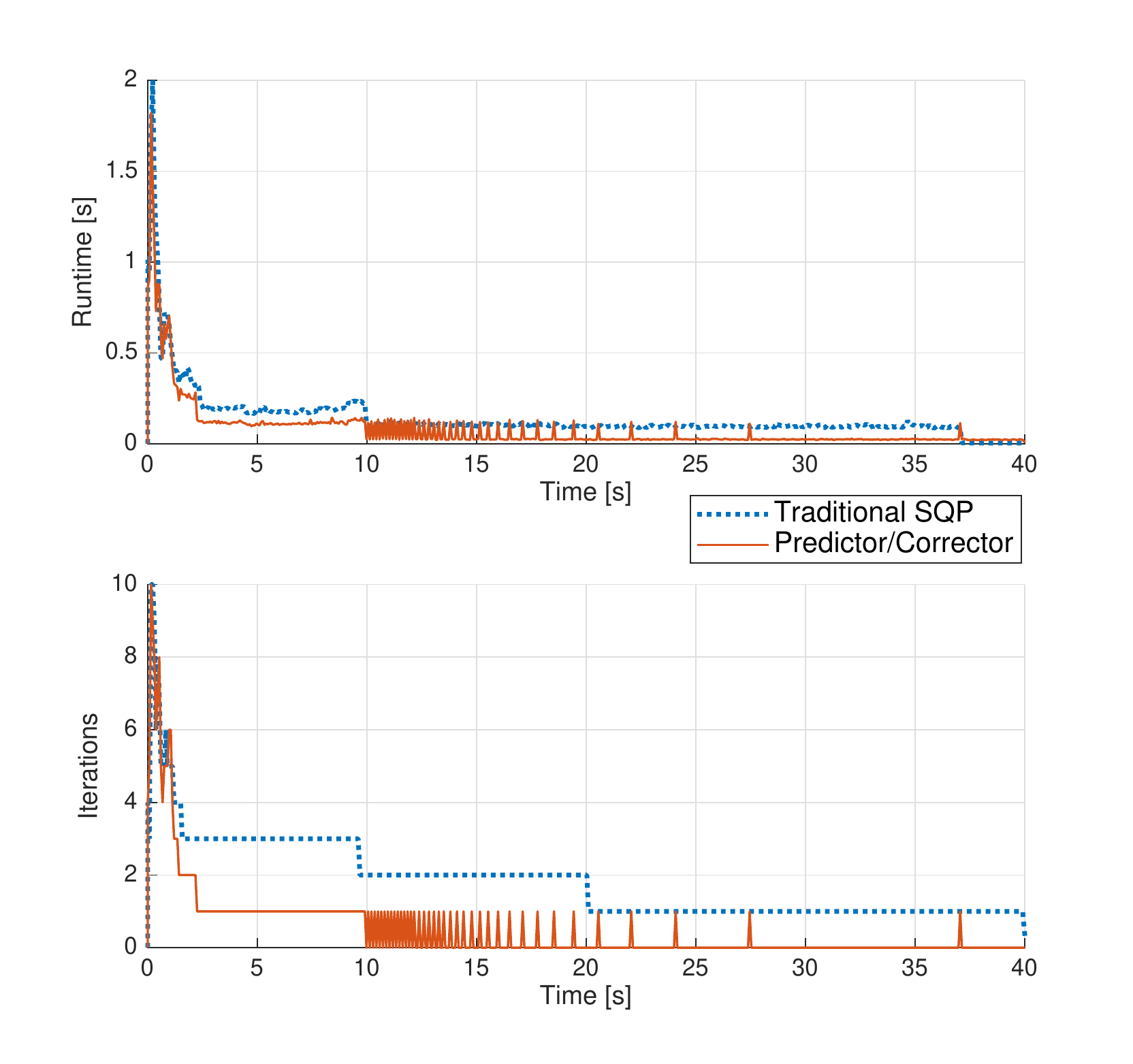}
  \caption{Runtime comparison for nonlinear MPC implemented using a one-step shift warmstart (dotted blue) and the proposed predictor-corrector (solid orange).}
  \label{fig:compare}
\end{figure}

\begin{figure}
  \centering
  \includegraphics[width=0.95\columnwidth]{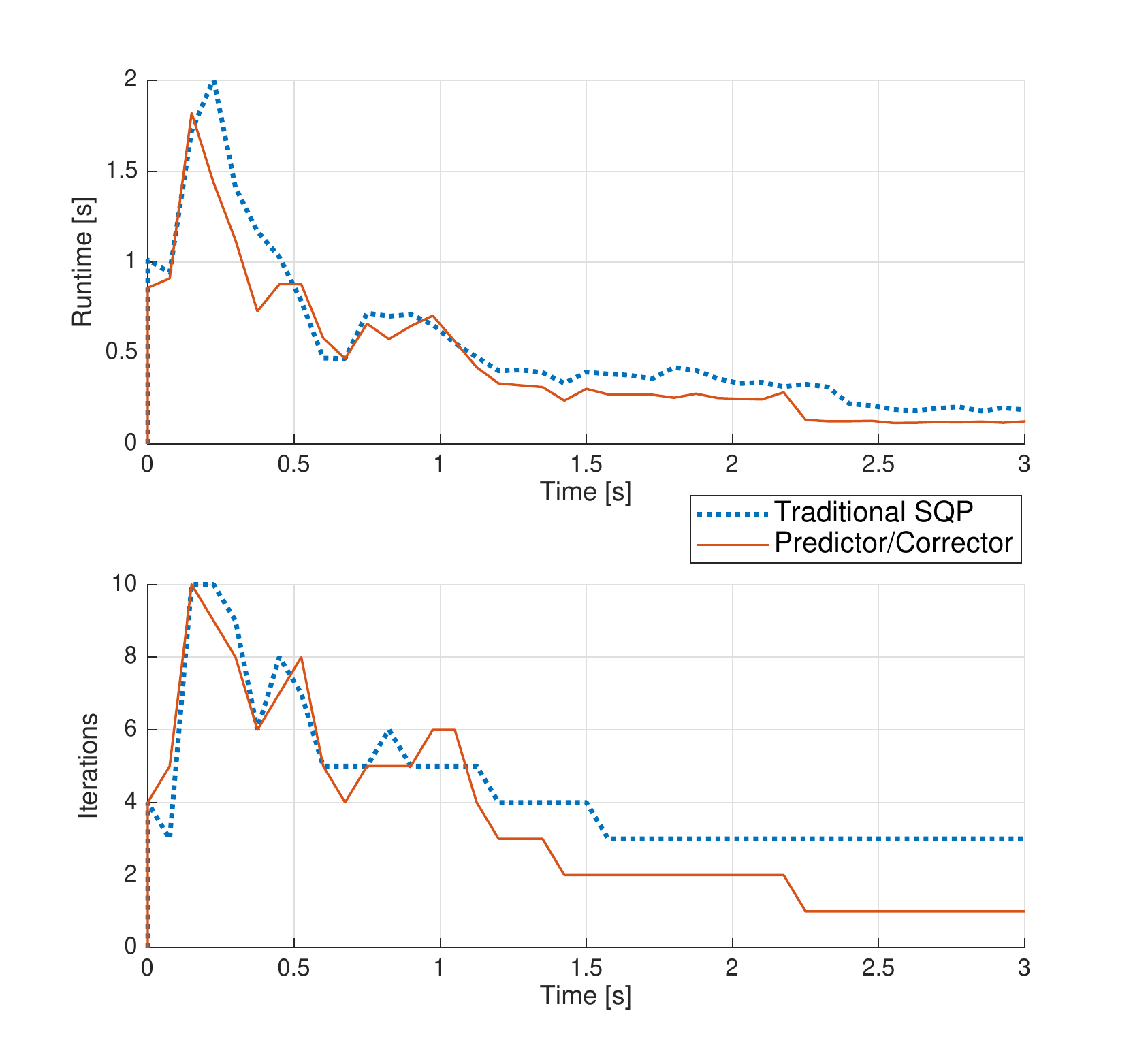}
  \caption{An enlarged section of a runtime comparison for nonlinear MPC implemented using a one-step shift warmstart (dotted blue) and the proposed predictor-corrector (solid orange).}
  \label{fig:compare2}
\end{figure}

Figure \ref{fig:response} displays the closed-loop response obtained using the MPC, which successfully enforces all the constraints. Note that the UAV position and velocity is reported for a longer horizon to show that the system is asymptotically stable, whereas all other signals are focused on the initial transient. Note that the LICQ is violated at $t = 0.6 s$. Figure~\ref{fig:compare} compares the computational cost obtained\footnote{The simulations were performed on a DELL Latitude 7390 2-in-1, Intel Core i7-8650U, 2.11 GHz, 16 GB laptop running MATLAB 2018b. All QPs were solved using MATLABs built in \texttt{quadprog} command, which implements an interior point algorithm, using default settings.} using a standard one-step shift warmstarting strategy, see e.g., \cite[Section 5.1]{magni2009nonlinear}, and the proposed predictor, Figure~\ref{fig:compare2} is an enlarged view of the first three seconds. The natural residual was subject to the tolerance value $\pi(p,z)\leq10^{-5}$.

The comparison in Figure \ref{fig:compare} shows that the predictor successfully achieves lower computational costs by reducing the number of SQP iterations performed by the solver. The average computational cost of one predictor step is in the order of $25\mu s$, compared to the computational cost of one SQP iteration which is around $100 \mu s$. This reduction in computational time is achieved without incurring significant drawbacks and without having to assume the LICQ, or indeed any CQ, for a system where that assumption is not guaranteed to hold. 

\section{Conclusions}
In this paper we proposed a new sensitivity-based warmstarting strategy for model predictive control of systems with nonlinear dynamics and linear state and control inequality constraints. The strategy involves a predictor which utilizes the semiderivative of the solution of the optimality  system. The method exploits the polyhedrality of the constraint set and requires fewer assumptions than comparable methods in the literature. Specifically, it does not require a difficult to verify constraint qualification. Numerical simulations demonstrate the potential of the strategy. Future work includes moving from polyhedral to more general constraints, relaxing the conditions for strong regularity, and developing tailored quadratic programming solvers for the predictor.

\bibliography{sens.bib}

\end{document}